\documentclass[11pt,twoside,english]{myamsart}
\advance\oddsidemargin by -1.0cm
\advance\evensidemargin by -1.0cm
\advance\topmargin by -1.0cm 
\usepackage{amssymb}
\usepackage{babel}
\usepackage{amstext}
\usepackage{amscd}   
\usepackage{epsfig}  
\usepackage{rotating}
\let\mathg\mathfrak
\theoremstyle{plain}

\newtheorem{lem}{Lemma}[section]
\newtheorem{thm}{Theorem}[section]            
\newtheorem{prop}{Proposition}[section]
\theoremstyle{definition}
\newtheorem{exa}{Example}[section]
\newtheorem{NB}{Remark}[section]

\renewcommand{\labelenumi}{$(\theenumi)$}
\newcommand{\bdm}{\begin{displaymath}}
\newcommand{\edm}{\end{displaymath}}
\newcommand{\be}{\begin{equation}}
\newcommand{\ee}{\end{equation}}
\newcommand{\ba}[1]{\begin{array}{#1}}
\newcommand{\ea}{\end{array}}

\newcommand{\btab}{\begin{tabular}}
\newcommand{\etab}{\end{tabular}}

\newcommand{\C}{\ensuremath{\mathbb{C}}}
\newcommand{\R}{\ensuremath{\mathbb{R}}}

\newcommand{\Q}{\ensuremath{\mathbb{Q}}}

\newcommand{\T}{\ensuremath{\mathrm{T}}}
\newcommand{\F}{\ensuremath{\mathrm{F}}}
\newcommand{\G}{\ensuremath{\mathrm{G}}}

\newcommand{\un}{\ensuremath{\mathg{u}}}
\newcommand{\su}{\ensuremath{\mathg{su}}}
\newcommand{\SU}{\ensuremath{\mathrm{SU}}}

\newcommand{\U}{\ensuremath{\mathrm{U}}}

\newcommand{\so}{\ensuremath{\mathg{so}}}
\newcommand{\hol}{\ensuremath{\mathg{hol}}}
\newcommand{\SO}{\ensuremath{\mathrm{SO}}}
\newcommand{\Spin}{\ensuremath{\mathrm{Spin}}}
\newcommand{\spin}{\ensuremath{\mathg{spin}}}

\newcommand{\g}{\ensuremath{\mathfrak{g}}}

\begin{document}
\def\haken{\mathbin{\hbox to 6pt{%
                 \vrule height0.4pt width5pt depth0pt
                 \kern-.4pt
                 \vrule height6pt width0.4pt depth0pt\hss}}}
    \let \hook\intprod
\setcounter{equation}{0}
%
%
\thispagestyle{empty}
%
\date{\today}
\title[$\G_2$-manifolds with parallel characteristic torsion]
{$\G_2$-manifolds with parallel characteristic torsion}
%
%
%
\author{Thomas Friedrich}
\address{\hspace{-5mm} 
{\normalfont\ttfamily friedric@mathematik.hu-berlin.de}\newline
Institut f\"ur Mathematik \newline
Humboldt-Universit\"at zu Berlin\newline
Sitz: WBC Adlershof\newline
D-10099 Berlin, Germany}
%
\thanks{Supported by the SFB 647 "Raum, Zeit, Materie" of the DFG}
\subjclass[2000]{Primary 53 C 25; Secondary 81 T 30}
\keywords{Cocalibrated $\G_2$-manifolds, 
connections with torsion}  
\begin{abstract}
We classify $7$-dimensional cocalibrated $\G_2$-manifolds with
parallel characteristic torsion and non-abelian holonomy. All these spaces
admit a metric connection $\nabla^{\mathrm{c}}$ with totally skew-symmetric
torsion and a spinor field $\Psi_1$ solving 
the equations in the common sector of type II superstring theory. 
There exist $\G_2$-structures with parallel characteristic
torsion that are not naturally reductive.
\end{abstract}
\maketitle
\pagestyle{headings}
%
%
%
%
%
%
\section{Metric connections with parallel torsion}\noindent
%
Consider a Riemannian manifold  and
denote by $\nabla^g$ its Levi-Civita connection.
Any $3$-form $\T$ defines via the formula
\bdm
\nabla_X Y \ := \ \nabla^g_X Y \, + \, \frac{1}{2} \, \T(X,Y, *)
\edm
a metric connection $\nabla$ with totally skew-symmetric torsion
$\T$. We are interested in the case that the torsion form is parallel,
$\nabla \T \, = \, 0$. 
Then $\T$ is coclosed, $\delta (\T) = 0$, and the differential $d \T$
depends only on the algebraic type of $\T$ (see \cite{FriedrichIvanov}),
\bdm
d \T \ = \ \sum_{i=1}^n (e_i \haken \T) \wedge (e_i \haken \T) \ .
\edm
The curvature tensor $\mathrm{R}^{\nabla} : \Lambda^2 \rightarrow
\Lambda^2$ of the connection $\nabla$ is symmetric,
$\mathrm{R}^{\nabla}(X,Y,U,V) = \mathrm{R}^{\nabla}(U,V,X,Y)$, see 
\cite{FriedrichIvanov}.
Moreover, if there exists a $\nabla$-parallel spinor field $\Psi$, 
we can compute the Ricci tensor directly using only the torsion form,
\bdm
2 \, \mathrm{Ric}^{\nabla}(X) \cdot \Psi \ = \ \big(X \haken d \T \big)
\cdot \Psi \, .
\edm
 In particular $\mathrm{Ric}^{\nabla}$
is parallel and divergence free,
$\nabla \mathrm{Ric}^{\nabla} = 0 \, , \
\mathrm{div}\big(\mathrm{Ric}^{\nabla}\big) = 0$ (see \cite{AFNP}). 
In \cite{AgFr2} we introduced a second order differential operator
$\Omega$ acting on spinors. It is a generalization of the Casimir operator
and its kernel contains all $\nabla$-parallel spinors. The formula
simplifies for parallel torsion, yielding
\bdm
\Omega \ = \ \Delta_{\T} \, + \, \frac{1}{16}( 2 \, \mathrm{Scal}^g \, + \, 
||\T||^2) \, - \, \frac{1}{4} \T^2  .
\edm
Consequently, any $\nabla$-parallel spinor $\Psi$ satisfies
\bdm
\T^2 \cdot \Psi \ = \ \frac{1}{4}( 2 \, \mathrm{Scal}^g \, + \, 
||\T||^2) \cdot \Psi \ .
\edm
In case we have more than  one $\nabla$-parallel spinor, the latter
equation
is an algebraic restriction for the $3$-form $\T$. Indeed, 
the endomorphism $\T^2$ acts as a scalar
on the space of all $\nabla$-parallel spinors.\\

\noindent
Almost Hermitian manifolds with parallel characteristic
torsion have been studied in \cite{AFS} and \cite{Sch}.
In this paper we consider $7$-dimensional Riemannian manifolds
equipped with a $\G_2$-structure $\varphi$. 
It is called \emph{cocalibrated} if the $3$-form
$\varphi$ satisfies the differential equation
\bdm
d\, * \, \varphi \ = \ 0  .
\edm
This ensures the existence of a  characteristic 
connection $\nabla^{\mathrm{c}}$ with totally skew-symmetric torsion
preserving the $\G_2$-structure  (see \cite{FriedrichIvanov})
\bdm
\T^{\mathrm{c}} \ = \ \frac{1}{6} \, ( d \varphi \, , \, * \varphi) 
\cdot \varphi \, - \, * \, d \varphi \,.
\edm
In general, the torsion form 
$\T^{\mathrm{c}}$ is not $\nabla^{\mathrm{c}}$-parallel. If it is,
there exists a spinor field $\Psi_1$ satisfying the
equations
\bdm
\nabla^{\mathrm{c}} \Psi_1 \ = \ 0 \, , \quad \T^{\mathrm{c}} \cdot \Psi_1 
\ = \ a \cdot \Psi_1
\, , \quad \delta(\T^{\mathrm{c}}) \ = \ 0 \ , \quad 
\nabla^{\mathrm{c}} \mathrm{Ric}^{\nabla^{\mathrm{c}}} \ = \ 0 \, ,
\edm
where the factor $a$ depends on the algebraic type of the torsion form.
Consequently, $\G_2$-manifolds with a parallel characteristic torsion
are solutions of the equations for the commen sector of type II superstring 
theory, see \cite{AFNP}.
The holonomy algebra $\hol(\nabla^{\mathrm{c}}) \subset 
\g_2$ preserves a $3$-form $ \T^{\mathrm{c}}$. The aim
of the present paper is the construction and classification
of cocalibrated $\G_2$-manifolds with parallel
characteristic torsion and non-abelian holonomy algebra. 
There are eight non-abelian subalgebras
of $\g_2$, see \cite{Dynkin}. 
For any of these algebras, we describe the set of admissible
torsion forms $\T^{\mathrm{c}}$. We discuss the  geometry
of the space $M^7$ in dependence of the type of its holonomy algebra
as well as the $\G_2$-orbit of the torsion form 
$ \T^{\mathrm{c}} \in \Lambda^3_1
\oplus \Lambda^3_{27}$. In particular, two subalgebras of $\g_2$ are realized
as holonomy algebras of a unique  cocalibrated 
$\G_2$-manifold with parallel characteristic
torsion, see Theorem \ref{irrso(3)} and 
Theorem \ref{Nonexistence}.
All other non-abelian subalgebras occur as holonomy
algebras for whole families of cocalibrated $\G_2$-manifolds. 

\textbf{Thanks.} I thank Ilka Agricola and Richard Cleyton (both from
Humboldt University Berlin) for many valuable discussions on the topic of this
article.

\section{The exceptional Lie algebra $\g_2$ and its subalgebras}\noindent
%
The group $\G_2$ is a compact, simple and simply-connected $14$-dimensional
Lie group. It consists of all elements in the group $\SO(7)$ preserving
the 3-form in seven dimensions
\begin{eqnarray*}
\varphi \ = \ e_{127} \, + \, e_{135} \, - \,  
e_{146} \, - \, e_{236} \, - \, e_{245} \, + \, e_{347} \, + \, e_{567} \, .
\end{eqnarray*}
The group $\G_2$ lifts into the spin group $\Spin(7)$ and, consequently, it
acts on spinors. Let $e_1,\ldots ,e_7$ be the standard orthonormal basis of 
Euclidian  space ${\Bbb R}^7$ and denote by $\mathrm{Cliff}({\Bbb R}^7)$
the real Clifford algebra. We will use  the following real representation of
this algebra on the space of real spinors $\Delta_7 := {\Bbb R}^8$:
\begin{eqnarray*}
e_{1} &=& \ \ E_{18} + E_{27} - E_{36} - E_{45} \, , \quad 
e_{2} \ = \ - E_{17} + E_{28} + E_{35} - E_{46} \, ,\\
e_{3} &=& - E_{16} + E_{25} - E_{38} + E_{47} \, , \quad
e_{4} \ = \ - E_{15} - E_{26} - E_{37} - E_{48} \, , \\
e_{5} &=& - E_{13} - E_{24} + E_{57} + E_{68} \, , \quad
e_{6} \ = \ \ \ E_{14} - E_{23} - E_{58} + E_{67} \, , \\
e_{7} &=& \ \ E_{12} - E_{34} - E_{56} + E_{78} \,
\end{eqnarray*}
where $E_{ij}$ denotes the standard basis of the Lie algebra $\so(8)$. We 
fix an orthonormal basis  $\psi _1 := (1,0,\ldots,0)\, , \ldots \, , 
\psi _8 := (0,0,\ldots,1)$ of spinors. The group $\G_2$ coincides with 
the subgroup of $\Spin(7)$ preserving the spinor  $\psi _1$.
Its Lie algebra ${\g_2}$ is the subalgebra of $\spin(7)$ containing all
forms $\omega = \sum \omega_{ij} \, e_{ij}$  such that the Clifford 
product $\omega \cdot \psi _1 = 0$
vanishes. This condition holds if and only if 
\begin{eqnarray*}
 \omega_{12} + \omega_{34} + \omega_{56} =0, & \quad -\omega_{13} +\omega_{24}
 -\omega_{67} =0 \, , & \quad -\omega_{14} - \omega_{23} - \omega_{57} =0 \, ,
 \\ 
-\omega_{16} - \omega_{25} + \omega_{37} =0 \, , & \quad \omega_{15}
 -\omega_{26} -\omega_{47} =0 \, , & \quad \omega_{17} + \omega_{36} + 
\omega_{45} =0 \, , \\ 
& \omega_{27} + \omega_{35} - \omega_{46} =0 \ . &
\end{eqnarray*}
These equations define the Lie algebra $\g_2$ inside $\spin(7)$. The
representations of $\g_2$ on $\R^7$ and on $\mathrm{Lin}(\psi_2, \psi_3,
\ldots , \psi_8)$ are equivalent.\\

\noindent
Dynkin's classical result that $\g_2$ has exactly
three maximal subalgebras of dimensions $8$, $6$ and $3$ respectively
implies that there are (up to conjugation) eight different non-abelian 
subalgebras of $\g_2$ (see \cite{Dynkin}). The $\G_2$-representation 
$\Lambda^3(\R^7)$ splits into a $1$-dimensional
trivial representation $\Lambda^3_1$ generated by the form $\varphi$, a
$7$-dimensional representation  $\Lambda^3_7$ containing all inner products
$e_i \haken * \varphi$, and the $27$-dimensional complement
$\Lambda^3_{27}$. Given a Lie subalgebra $\g \subset 
\g_2$ let $(\Delta_7)_{\g}$ and $(\Lambda^3_{27})_{\g}$ be the
space of all $\g$-invariant spinors and the space of all $\g$-invariant
$3$-forms in $\Lambda^3_{27}$, respectively. The space of torsion
forms of cocalibrated $\G_2$-manifolds  
with parallel torsion is the following set,
\bdm
 \mathrm{Tor}_{\g} \ := \, \big\{ \mathrm{T} \in \Lambda^3_1 \oplus 
(\Lambda^3_{27})_{\g}\, :
   \, \mathrm{T}^2 \ \mathrm{acts} \ \mathrm{as} \ \mathrm{a} 
\ \mathrm{scalar} 
\ \mathrm{on} \ (\Delta_7)_{\g} \big\} \, . 
\edm
Since the $\G_2$-orbit of the characteristic torsion is a geometric
invariant, two torsion forms $\mathrm{T}_1$ and $\mathrm{T}_2$ define 
equivalent geometric structures if they are equivalent as 
$3$-forms under the action of the
normalizer $\mathcal{N}(\G)$ of the Lie group $\G \subset \G_2$. Consequently,
the relevant set is $ \mathrm{Tor}_{\g}/\mathcal{N}(\G)$.
We computed these spaces for any of the eight non-abelian subalgebras
of $\g_2$. Here we formulate 
only the result of these computations.
\vspace{3mm}

\noindent
{\bf The subalgebra $\su(3) \subset \g_2$.}\noindent
%
\vspace{2mm}

\noindent
The Lie algebra preserving two spinors $\psi_1$ and $\psi_2$ is isomorphic
to $\su(3)$. It is the unique subalgebra of $\g_2$ of dimension eight.
The representation on Euclidian space splits 
into $\R^7 = \C^3 \oplus \R^1$. $(\Lambda^3_{27})_{\su(3)}$ is
one-dimensional and generated by
\bdm
\T \ = \ 4 \, e_{127} \, - \, 3 \, e_{135} \, + \, 3 \, e_{146} \, + \, 
3 \, _{236} \, + \, 3 \, e_{245} \, + \, 4 \, e_{347} \, + \, 4 \, e_{567}.
\edm
The set $\mathrm{Tor}_{\su(3)}$ is the union of two lines, 
\bdm
\mathrm{Tor}_{\su(3)}  \ = \ \big\{a \cdot (e_{127} \, + \, e_{347} \, + 
\, e_{567}) \big\} \ \cup \ \big\{b \cdot (- \, e_{135} \, + \, e_{146} \, + 
\, e_{236} \, + \, e_{245} \big\})  .
\edm
Remark that the intersections  $\mathrm{Tor}_{\su(3)} \cap \Lambda^3_1 = 0$
and  $\mathrm{Tor}_{\su(3)} \cap \Lambda^3_{27} = 0$ are trivial.
\vspace{3mm}

\noindent
{\bf The subalgebra $\un(2) \subset \su(3) \subset \g_2$.}\noindent
%
\vspace{2mm}

\noindent
This Lie subalgebra is generated by 
\bdm
P_1 \ := \ e_{13} \, + \, e_{24} \, , \quad
P_2 \ := \ e_{14} \, - \, e_{23} \, , \quad
P_3 \ := \ e_{12} \, - \, e_{34} \, , \quad
Q_3 \ := \ e_{12} \, + \, e_{34} \, - \, 2 e_{56} \ .
\edm
The subalgebra $\un(2) \subset \g_2$ preserves two spinors $\psi_1 \, , \, 
\psi_2$ and it acts on Euclidian space as  
$\R^7 = \C^2 \oplus \C \oplus \R^1$. The space
$(\Lambda^3_{27})_{\un(2)}$ is two-dimensional and parameterized by
\bdm
\T_{a,b} \ := \ (a \, + \, b) \, (- \, e_{135} \, + \, e_{146} \, + \, 
e_{236} \, + \, e_{245}) \, + \, 2a \, (e_{347} \, + \, e_{127}) \,
+ \, 4b \, e_{567} .
\edm
$\mathrm{Tor}_{\un(2)}$ is the union of two planes 
in $\Lambda^3_1 \oplus (\Lambda^3_{27})_{\un(2)}$. The first plane $\Pi_1$
is parameterized by
\bdm
\T_{a,b} \, + \, c \cdot \varphi \, , \quad \mathrm{where} \quad a \, + \, b
\ = \ c \ . 
\edm
The second plane $\Pi_2$ is the family
\bdm
\T_{a,b} \, + \, c \cdot \varphi \, , \quad \mathrm{where} \quad 4\, (a \, + \, b)
\ = \ - \, 3 \, c \ .
\edm
Here the intersection $\Pi_1 \cap \Pi_2$ is
a line in $(\Lambda^3_{27})_{\un(2)}$.
\vspace{3mm}

\noindent
{\bf The subalgebra $\su(2) \subset \su(3) \subset \g_2$.}\noindent
%
\vspace{2mm}

\noindent
This Lie subalgebra is generated by $P_1, P_2, P_3$ and stabilizes
 four spinors $\psi_1, \psi_2,
\psi_3, \psi_4$. Under $\su(2)$, we have the decomposition
$\R^7 = \C^2 \oplus \R^1 \oplus \R^1 \oplus \R^1$. It turns out that 
$(\Lambda^3_{27})_{\su(2)}$ is a $6$-dimensional space. Moreover, the 
normalizer $\mathcal{N}(\SU(2))$ is a three-dimensional subgroup. It acts
on the spinor $\psi_2 , \psi_3, \psi_4$ in a non-trivial way.
Consequently, the geometrically
relevant set $\mathrm{Tor}_{\su(2)}/\mathcal{N}(\SU(2))$ consists of
$3$-forms admitting the four spinors $\psi_1 , \psi_2 , \psi_3 , 
\psi_4$ as eigenspinors. It is the union of eight lines in 
$\Lambda^3_1 \oplus (\Lambda^3_{27})_{\su(2)}$ generated
by the following forms:
\begin{eqnarray*}
\T &=& e_{127} \, - \, e_{135} \, + \, e_{146} \, + \, e_{236} \, + \, e_{245} \, + \, e_{347} \, - \,2 \, e_{567} \ , \\
\T &=& e_{127} \, + \, e_{135} \, - \, e_{146} \, - \, e_{236} \, - \, e_{245}
\, + \, e_{347} \, - \,2 \, e_{567} \ , \\
\T &=& e_{127} \, + \, e_{135} \, + \, e_{146} \, + \, e_{236} \, - \, e_{245}
\, + \, e_{347} \, + \,2 \, e_{567} \ , \\
\T &=& e_{127} \, - \, e_{135} \, - \, e_{146} \, - \, e_{236} \, + \, e_{245}
\, + \, e_{347} \, + \,2 \, e_{567} \ , \\
\T &=& e_{135} \, - \, e_{245} \, , \quad \T \ = \ e_{146} \, + \, e_{236} \,
, \quad \T \ = \ e_{127} \, + \, e_{347} \, , \quad \T \ = \ e_{567}  .
\end{eqnarray*}
The intersections of  $\mathrm{Tor}_{\su(2)} \cap \Lambda^3_1 = 0$
and  $\mathrm{Tor}_{\su(2)} \cap \Lambda^3_{27} = 0$ are trivial.
\vspace{3mm}

\noindent
{\bf The subalgebra $\su_c(2) \subset \g_2$.}\noindent
%
\vspace{2mm}

\noindent
The centralizer of the subalgebra $\su(2) \subset \g_2$ is a subalgebra
of $\g_2$ which is isomorphic, but not conjugated to $\su(2)$. We denote
this algebra by  $\su_c(2)$. It is generated by 
\bdm
Q_1 \ := \ - \, e_{14} \, - \, e_{23} \, + \, 2 e_{57} \, , \quad
Q_2 \ := \ - \, e_{13} \, + \, e_{24} \,  + \, 2 e_{67} \, , \quad
Q_3 \ := \ e_{12} \, + \, e_{34} \, - \, 2 e_{56}  .
\edm 
The subalgebra $\su_c(2) \subset \g_2$ preserves only one spinor 
$\psi_1$ and Euclidian space splits under its action into 
$\R^7 = \C^2 \oplus \R^3$. The space
$(\Lambda^3_{27})_{\su_c(2)}$ is one-dimensional and generated by the
$3$-form
\bdm
\T \ = \ \varphi \, - \, 7 \, e_{567}  .
\edm
The set  $\mathrm{Tor}_{\su_c(2)}$ coincides with $\Lambda^3_1 \oplus
(\Lambda^3_{27})_{\su_c(2)}$ and is generated by $\varphi$ and $e_{567}$.
\vspace{3mm}

\noindent
%
{\bf The subalgebra $\R^1 \oplus \su_c(2) \subset \g_2$.}\noindent
%
\vspace{2mm}

\noindent
This subalgebra is generated by $P_1, Q_1, Q_2, Q_3$. 
The subalgebra $\R^1 \oplus \su_c(2) \subset \g_2$ preserves only one spinor 
$\psi_1$, and we have $\R^7 = \C^2 \oplus \R^3$. The space
$(\Lambda^3_{27})_{\R^1 \oplus \su_c(2)}$ is one-dimensional and 
generated by the $3$-form
\bdm
\T \ = \ \varphi \, - \, 7 \, e_{567}  .
\edm
The set  $\mathrm{Tor}_{\R^1 \oplus \su_c(2)}$ is generated by 
$\varphi$ and $e_{567}$.
\vspace{3mm}

\noindent
{\bf The subalgebra $\su(2) \oplus \su_c(2) \subset \g_2$.}\noindent
%
\vspace{2mm}

\noindent
The subalgebra $\su(2) \oplus \su_c(2)$ is generated by 
$P_1, P_2, P_3, Q_1, Q_2, Q_3$ and preserves the spinor $\psi_1$,
the representation on $\R^7$ decomposes into $\R^7 = \R^4 \oplus \R^3$.
The space
$(\Lambda^3_{27})_{\su(2) \oplus \su_c(2)}$ is one-dimensional and 
generated by the $3$-form
\bdm
\T \ = \ \varphi \, - \, 7 \, e_{567} , 
\edm
and the set  $\mathrm{Tor}_{\su(2) \oplus \su_c(2)}$ is generated by 
$\varphi$ and $e_{567}$.
\vspace{3mm}

\noindent
{\bf The subalgebra $\so(3) \subset \su(3) \subset \g_2$.}\noindent
%
\vspace{2mm}

\noindent
This subalgebra is generated by 
\bdm
S_1 \ := \ e_{12} \, - \, e_{56}  , \quad
S_2 \ := \ e_{13} \, + \, e_{24} \, + \, e_{35} \, + \, e_{46}  , \quad
S_3 \ := \ e_{14} \, - \, e_{23} \, + \, e_{36} \, - \, e_{45}  
\edm
and stabilizes $\psi_1$ and $\psi_2$. We have the further splitting
$\R^7 = \R^3 \oplus \R^3 \oplus \R^1$, and  $(\Lambda^3_{27})_{\so(3)}$ 
is spanned by three linearly independet $3$-forms as follows,
\begin{eqnarray*}
\T_{a,b,c} \ := && a\, (- 2 \, e_{123} \, + \, e_{136} \, - \, e_{145} \, + \,
e_{235} \, + \, e_{246} \, + \, 2 \, e_{356}) \\ 
&+& \, b \, (- 2 \, e_{124} 
\, - \, e_{135} \, - \, e_{146} \, + \,
e_{236} \, - \, e_{245} \, + \, 2 \, e_{456}) \\ 
&+& \, c \, (4 \, e_{127} 
\, - \, 3 \, e_{135} \, + \, 3\, e_{146} \, + \,
3 \, e_{236} \, + \, 3 \, e_{245} \, + \, 4 \, e_{347} \, + \, 4 \, e_{567})
\ .
\end{eqnarray*} 
$\mathrm{Tor}_{\so(3)}$ is the union of two hyperplanes in 
$\Lambda^3_1 \oplus (\Lambda^3_{27})_{\so(3)}$. The intersection
$\mathrm{Tor}_{\so(3)} \cap \Lambda^3_{27}$ is two-dimensional and consists
of all forms with $c=0$, $\T_{a,b,0}$.

\vspace{3mm}

\noindent
{\bf The subalgebra $\so_{ir}(3) \subset \g_2$.}\noindent
%
\vspace{2mm}

\noindent
The unique irreducible $7$-dimensional
real representation of $\so(3)$ is contained
in $\g_2$. Consequently we obtain a maximal subalgebra of dimension three,
$\so(3)_{ir} \subset \g_2$. 
This algebra preserves only one spinor. 
$(\Lambda^3_{27})_{\so_{ir}(3)}= 0$ is trivial
and $\mathrm{Tor}_{\so_{ir}(3)} = \Lambda^3_1$ is generated by 
$\varphi$.
%
%
\section{The characteristic connection of a cocalibrated $\G_2$-manifold}\noindent
A $\G_2$-manifold $(M^7,g,\varphi)$ is called cocalibrated if the $3$-form
$\varphi$ satisfies the differential equation
\bdm
d\, * \, \varphi \ = \ 0  ,
\edm
see \cite{FerGray}.
Then there exists a unique connection $\nabla^{\mathrm{c}}$ preserving the
$\G_2$-structure with totally skew-symmetric torsion (see 
\cite{FriedrichIvanov})
\bdm
\T^{\mathrm{c}} \ = \ \frac{1}{6} \, ( d \varphi \, , \, * \varphi) \cdot
\varphi \, - \, * \, d \varphi.
\edm
The torsion form belongs to $\Lambda^3_1 \oplus
\Lambda^3_{27}$. If the torsion is parallel, 
the Ricci tensor depends only on the
algebraic type of the torsion form  $\T^{\mathrm{c}}$, see
\cite{FriedrichIvanov}. In particular, $\mathrm{Ric}^{\nabla^{\mathrm{c}}}$
is parallel and divergence free, see \cite{AFNP},
\bdm
\nabla^{\mathrm{c}} \mathrm{Ric}^{\nabla^{\mathrm{c}}} \ = \ 0 \, , \quad
\mathrm{div}\big(\mathrm{Ric}^{\nabla^{\mathrm{c}}}\big) \ = \ 0 \ . 
\edm
On all these spaces, there exists a spinor field $\Psi_1$ satisfying the
equations
\bdm
\nabla^{\mathrm{c}} \Psi_1 \ = \ 0  , \quad \T^{\mathrm{c}} \cdot \Psi_1 
\ = \ a \cdot \Psi_1, \quad \delta(\T^{\mathrm{c}}) \ = \ 0  , 
\edm
where the factor $a$ depends on the algebraic type of the torsion form.\\

\noindent
The subclass of nearly parallel $\G_2$-structures is given by
the condition that the $4$-forms $d \varphi$ and $* \varphi$ are proportional.
In this case the characteristic torsion $\T^{\mathrm{c}}$ is proportional
to $\varphi$ and, consequently, it is automatically parallel, 
$\nabla^{\mathrm{c}}
\T^{\mathrm{c}} = 0$. Although $\su(3)$ preserves the spinors $\psi_1$ and
$\psi_2$,  $\varphi$ acts on these spinors with different eigenvalues, hence
no  subalgebra of $\su(3)$ can  occur as holonomy algebra 
$\hol(\nabla^{\mathrm{c}})$ of the characteristic connection.
\vspace{2mm}

\begin{prop} \label{Nearlyparallel}
Let $(M^7,g,\varphi)$ be a nearly parallel $\G_2$-manifold which is not
parallel, $\T^{\mathrm{c}} \neq 0$. Then the holonomy algebra
$\hol(\nabla^{\mathrm{c}})$ is not a subalgebra of $\su(3)$. The
$\nabla^{\mathrm{c}}$-parallel spinor of the nearly parallel 
$\G_2$-manifold is unique.
\end{prop}
\vspace{2mm}

\noindent
On the other hand, nearly parallel $\G_2$-manifolds with characteristic
holonomy $\hol(\nabla^{\mathrm{c}}) = \su(2) \oplus \su_c(2)$ exist.
Indeed, consider a $3$-Sasakian manifold. By an appropriate rescaling of
its metric in the direction of the three-dimensional bundle
spanned by the contact forms, one obtains a nearly parallel
$\G_2$-manifold of that type (see \cite{FKMS}, Theorem 5.4). The irreducible
naturally reductive homogeneous space $\SO(5)/\SO_{ir}(3)$ is an example 
of a nearly parallel $\G_2$-manifold with characteristic holonomy  
$\hol(\nabla^{\mathrm{c}}) = \so_{ir}(3)$.
\vspace{2mm}

\noindent
We now  discuss $\G_2$-manifolds with  characteristic torsion in
$\Lambda^3_{27}$ (structures of pure type
$\mathcal{W}_3$, see \cite{FerGray}). 
The differential equations characterizing these structures
are
\bdm
d\, * \, \varphi \ = \ 0  , \quad   ( d \varphi \, , \, * \varphi) \ = 
\ 0  . 
\edm
In this case the formula for the characteristic torsion simplifies,
\bdm
\T^{\mathrm{c}} \ = \ - \, * \, d \varphi  .
\edm
Observe that a cocalibrated $\G_2$-manifold is of type $\mathcal{W}_3$ if and
only if the Clifford product $\T^{\mathrm{c}} \cdot \Psi = 0$ of the
characteristic torsion and the canonical $\nabla^{\mathrm{c}}$-parallel spinor
vanishes.
Since $\T^{\mathrm{c}} \neq 0$ is preserved by the characteristic
connection, the holonomy algebra  $\hol(\nabla^{\mathrm{c}})$ is a proper
subalgebra of $\g_2$. Moreover, the intersection  
 $\mathrm{Tor}_{\g} \cap \Lambda^3_{27} = 0$ is trivial for $\g = \su(3),
 \su(2)$ and $\so_{ir}(3)$. The algebraic computations 
yield now the following result.
\vspace{2mm}

\begin{prop}
Let $M^7,g,\varphi)$ be a $\G_2$-manifold of type $\mathcal{W}_3$ with 
parallel characteristic torsion. If  $\hol(\nabla^{\mathrm{c}}) = \su(3),
\su(2)$ or $\so_{ir}(3)$, then $M^7$ is a parallel $\G_2$-manifold, 
$\T^{\mathrm{c}} = 0$.
\end{prop}
\vspace{2mm}

\noindent
We will study the cocalibrated $\G_2$-geometries with  non-abelian holonomy 
algebra $\hol(\nabla^{\mathrm{c}})$ case by case.
%
\section{$\G_2$-manifolds with parallel torsion and 
$\hol(\nabla^{\mathrm{c}}) = 
\su(3)$}\noindent
A cocalibrated $\G_2$-manifolds with non-trivial parallel torsion and 
$\hol(\nabla^{\mathrm{c}}) = \su(3)$ admits a $\nabla^{\mathrm{c}}$-parallel
vector field $e_7$ and two $\nabla^{\mathrm{c}}$-parallel real
spinor fields $\Psi_1 , \ \Psi_2 = e_7 \cdot \Psi_1$. Let us introduce 
the following two globally well defined and  
$\nabla^{\mathrm{c}}$-parallel forms,
\bdm
\F \ := \ e_{12} \, + \, e_{34} \, + \, e_{56}  , \quad
\Sigma \ := \ e_{135} \, - \, e_{146} \, - \, e_{236} \, - \, e_{245} 
\edm
and define a  $\G_2$-structure  by setting
\bdm
\varphi \ = \ \Sigma \, + \,  \F \wedge e_7. 
\edm
The torsion form
$\T^{\mathrm{c}}$ acts on these spinor fields with the same or with
opposite eigenvalues. First we discuss the case that
the torsion form acts on both spinor fields with the same eigenvalue. 
Up to a scaling of the metric we can assume that 
$\T^{\mathrm{c}}$ is given by the formula
\bdm
\T^{\mathrm{c}} \ = \ 2 \, (e_{127} \, + \, e_{347} \, +  \, e_{567}) \ = \ 2
\, \F \wedge e_7 \, ,
\quad \T^{\mathrm{c}} \cdot \Psi_1 \ = \ - \, 6 \, \Psi_1 \, , \quad
\T^{\mathrm{c}} \cdot \Psi_2 \ = \ - \, 6 \, \Psi_2 \ .
\edm  
The equation $\nabla^{\mathrm{c}} e_7 = 0$ yields directly
\bdm
\nabla^{g}_X e_7 \ = \ X \haken \F \, , \quad d e_7 \ = \ 2 \, \F \, , \quad
\T^{\mathrm{c}} \ = \ d e_7 \wedge e_7 \ . 
\edm
Consequently, the tuple $(M^7, g , e_7, \F)$ is a Sasakian manifold and
$\T^{\mathrm{c}}$ is the characteristic torsion of this contact structure, 
see \cite{FriedrichIvanov}. Consider the  $\nabla^{\mathrm{c}}$-parallel
spinors $\Psi_1, \, \Psi_2$. Using the special algebraic 
formula for the torsion we compute that both spinors are solutions of the equation
\bdm
\nabla^g_X \Psi \ = \ - \, \frac{1}{2} \, X \cdot \Psi \, - \, g(X,e_7) \cdot
e_7 \cdot \Psi \ . 
\edm
In the paper \cite{FKim} we discussed the integrability conditions for this 
system. In the notations of this paper, the spinors $\Psi_1, \Psi_2$ are 
Sasakian quasi-Killing spinors of type $( -  1/2 , - 1)$.
Then $(M^7, g , e_7, \F)$ has to be a $\eta$-Einstein Sasakian manifold with
Ricci tensor
\bdm
\mathrm{Ric}^g = \ 10 \cdot g \, - \, 4 \cdot e_7 \otimes e_7  .
\edm
The $3$-form $\Sigma$ can be interpreted in contact geometry, too. It
is a horizontal,  $\nabla^{\mathrm{c}}$-parallel and belongs
to the $\Lambda^3_2(\R^6)$-component in the
decomposition of $\Lambda^3(\R^6)$ under the action of the group $\U(3)$,
see \cite{AFS}. Let us fix now a  simply-connected 
$\eta$-Einstein Sasakian manifold of the prescribed type. Then there exist
two Sasakian quasi-Killing spinors of type $(-1/2 ,- 1)$ 
(see \cite{FKim}, Theorem 6.3). 
Moreover, we can reconstruct the $\G_2$-structure as well as the
 $\nabla^{\mathrm{c}}$-parallel spinor fields. Finally we obtain the following
 result.

\begin{thm}
Any cocalibrated $\G_2$-manifold such that
the characteristic torsion acts on both $\nabla^{\mathrm{c}}$-parallel
spinors by the same eigenvalue and
\bdm
\nabla^{\mathrm{c}} \T^{\mathrm{c}} \ = \ 0 \, , \quad 
\T^{\mathrm{c}} \ \neq \ 0 \, , \quad
\hol(\nabla^{\mathrm{c}}) \ = \ \su(3) \
\edm
holds is homothetic to an  $\eta$-Einstein Sasakian manifold. Its
Ricci tensor is given by the formula
\bdm
\mathrm{Ric}^g = \ 10 \cdot g \, - \, 4 \cdot e_7 \otimes e_7 \ .
\edm
Conversely, a simply-connected  $\eta$-Einstein Sasakian manifold with Ricci
tensor $\mathrm{Ric}^g =  10 \cdot g \, - \, 4 \cdot e_7 \otimes e_7$ 
admits a cocalibrated $\G_2$-structure with parallel characteristic torsion
and characteristic holonomy contained in $\su(3)$.
\end{thm}
\vspace{2mm}

\noindent
The paper \cite{BGM} is an introduction to the geometry of 
$\eta$-Einstein Sasakian manifolds.
\vspace{3mm}

\noindent
Next we investigate the case that
the torsion form acts on the parallel spinor fields with opposite
eigenvalues. 
Then we have
\bdm
\T^{\mathrm{c}} \ = \ e_{135} \, - \, e_{146} \, -  \, e_{236} \, - \, 
e_{245} \ = \ \Sigma \, ,
\quad \T^{\mathrm{c}} \cdot \Psi_1 \ = \ - \, 4 \, \Psi_1 \, , \quad
\T^{\mathrm{c}} \cdot \Psi_2 \ = \ 4 \, \Psi_2 \ ,
\edm  
and the scalar curvature is $\mathrm{Scal}^g = 30$. Since
$\T^{\mathrm{c}}$
does not depend on $e_7$ and $\nabla^{\mathrm{c}} e_7 = 0$, we conclude that
$e_7$ is parallel with respect to the Levi-Civita connection,
$\nabla^{g} e_7 = 0$. The manifold $M^7$ splits isometrically
into $M^7 = X^6 \times \R^1$, where $X^6$ is an almost complex
manifold with K\"ahler form $\F$ and characteristic torsion $\T^{\mathrm{c}}$.
This torsion form belongs to the $\Lambda^3_2(\R^6)$-component in the
decomposition of $\Lambda^3(\R^6)$ under the action of the group $\U(3)$,
see \cite{AFS}. Consequently, $X^6$ is a nearly K\"ahler
(non-K\"ahlerian) $6$-manifold. Vice versa, we can reconstruct the
$7$-dimensional $\G_2$-structure out of the $6$-dimensional nearly
K\"ahler structure by the formula
\bdm
\varphi \ = \ \Sigma  +  \F \wedge e_7  ,
\edm
where $\F$ is the K\"ahler form and $\Sigma$ the parallel
$3$-form of type $\Lambda^3_2$ of $X^6$. Let us summarize the result.
\begin{thm}
A complete, simply-connected cocalibrated $\G_2$-manifold such that
the characteristic torsion acts on  $\nabla^{\mathrm{c}}$-parallel
spinors by opposite eigenvalues and
\bdm
\nabla^{\mathrm{c}} \T^{\mathrm{c}} \ = \ 0  , \quad 
\T^{\mathrm{c}} \ \neq \ 0  , \quad
\hol(\nabla^{\mathrm{c}}) \ = \ \su(3) 
\edm
holds is isometric the the product of a nearly K\"ahler $6$-manifold
by $\R$. Conversely, any such product admits a cocalibrated 
$\G_2$-structure with parallel torsion and holonomy contained in $\su(3)$.
\end{thm}
%
%
\section{$\G_2$-manifolds with parallel torsion and 
$\hol(\nabla^{\mathrm{c}}) = 
\un(2)$}\noindent
In this case the following forms are globally well defined and 
$\nabla^{\mathrm{c}}$-parallel,
\bdm
e_7 \, , \quad \Omega_1 \ : = \ e_{12} \, + \, e_{34} \, , \quad
\Omega_2 \ : = \ e_{56} \, , \quad \Sigma \ := \ e_{135} \, - \, 
e_{146} \, - \, e_{236} \, - \, e_{245} . 
\edm
The $\G_2$-structure as well as the admissible torsion forms are given
by the formulas
\begin{eqnarray*}
\varphi &=& \Sigma  +  \Omega_1 \wedge e_7 +  \Omega_2 \wedge e_7, \\
\T^{\mathrm{c}} &=& -  (a  +  b  -  c)  \Sigma  +  (2a +  c) \, 
\Omega_1 \wedge e_7 +  (4b  +  c)  \Omega_2 \wedge e_7  ,
\end{eqnarray*}
where one of the conditions $ a  +  b = c  $ or $4 (a  +  b) = - 3c$
is satisfied. Using Proposition 4.2 of the paper \cite{AgFr1} we compute
the differentials of these forms.
\begin{lem}
\begin{eqnarray*}
de_7 &=& (2 a \, + c) \, \Omega_1 \, + \, (4b \, + \, c) \, \Omega_2 \, , \\
d\Omega_1 &=& 2 \, (a \, + \, b \, - \, c) \, (e_7 \haken * \Sigma) \ = \
2 \, d \Omega_2 \, , \\
d \Sigma &=& - \, 4 \, (a \, + \, b \, - \, c) \, \Omega_1 \wedge \Omega_2 
\, + \, ( 4a \, + \, 4b \, + \, 3c) \, (* \Sigma) \, , \\
d* \Sigma &=& 0 \, , \quad  d(e_7 \haken * \Sigma) \ = \ 
 ( 4a \, + \, 4b \, + \, 3c) \, \Sigma \wedge e_7  .
\end{eqnarray*}
\end{lem} 

\noindent
First we study the case $ 4a  + 4b +  3c = 0, c  \neq  0$. 
The torsion form is given by
\bdm
\T^{\mathrm{c}} 
\ = \ \frac{7}{4} c \Sigma  +  (2 a  +  c)  (\Omega_1  - 
 2  \Omega_2) \wedge e_7
\edm 
and it acts on the spinors $\Psi_1$ and $\Psi_2$ with opposite
eigenvalues $\pm \, 7 c$. We compute the scalar curvature,
\bdm
\mathrm{Scal}^g \ = \ -  12 a^2  -  12 a c  +  \frac{711}{8} c^2  .
\edm
The formulas for the differentials of the globally defined forms can be
simplified,
\begin{eqnarray*}
de_7 &=& (2 a \, + c) \, (\Omega_1  -  2 \Omega_2)  , \quad
d\Omega_1 \ = \ -  \frac{7}{2}  c \, (e_7 \haken * \Sigma) \ = \ 2 
d \Omega_2 , \\
d* \Sigma &=& 0 \, , \quad  d(e_7 \haken * \Sigma) \ = \ 0 \, , \quad
d \Sigma \ = \ 7 c   \Omega_1 \wedge \Omega_2  .
\end{eqnarray*}
In particular, all possible Lie derivatives vanish
\bdm
\mathcal{L}_{e_7} \Omega_1 \ = \ 0 \, , \quad
\mathcal{L}_{e_7} \Omega_2 \ = \ 0 \, , \quad
\mathcal{L}_{e_7} \Sigma \ = \ 0 \, , \quad
\mathcal{L}_{e_7} (e_7 \haken * \Sigma) \ = \ 0 \, , \quad
\mathcal{L}_{e_7} \T^{\mathrm{c}} \ = \ 0  .
\edm
Let us discuss the regular case, i.\,e.~we assume that $e_7$ induces a free
action of the group $\mathrm{S}^1$. Then $\pi : M^7 \longrightarrow
\tilde{X}^6$ is a principal fiber bundle over a smooth manifold 
$\tilde{X}^6$. Moreover, on this manifold there exist differential forms
$\tilde{\Omega}_1 \, \tilde{\Omega}_2$ and  $\tilde{\Sigma}$ such that
\bdm
\tilde{\T}^{\mathrm{c}} \ = \ \frac{7}{4} \, c \, \tilde{\Sigma} \, , \quad
\tilde{*}\tilde{\Sigma} \ = \ - \, \tilde{e_7 \haken * \Sigma }
\edm
holds, where $\tilde{*}$ denotes the Hodge operator of $\tilde{X}^6$.
We introduce the form
\bdm
\tilde{\Omega} \ := \ \tilde{\Omega}_1  +  \tilde{\Omega}_2  .
\edm
Then we obtain
\bdm
d \tilde{\Omega} \ = \ 3 \frac{7}{4}  c  \tilde{*} \tilde{\Sigma}
\ = \ 3  \tilde{*} \tilde{\T}^{\mathrm{c}}  . 
\edm
The torsion form $\tilde{\T}^{\mathrm{c}} \neq 0$ 
is of type $\Lambda^3_2$ in the sense
of almost Hermitian geometry on $\tilde{X}^6$. The last equation shows
that  $\tilde{X}^6$ is a nearly K\"ahler manifold (see \cite{AFS}, section
4.2.) with reduced characteristic holonomy 
$\hol(\tilde{\nabla}^{\mathrm{c}}) = \un(2) \subset \su(3)$. Then 
$\tilde{X}^6$ is isomorphic to the projective space $\mathbb{CP}^3$ or to
the flag manifold $\mathbb{F}(1,2)$ equipped with their standard nearly
K\"ahler structure coming from the twistor construction, see \cite{BM}.
The form $\tilde{\Omega}_1  -  2  \tilde{\Omega}_2$ is their
standard K\"ahler form. Conversely, if $\tilde{X}^6 = \mathbb{CP}^3 ,  
\mathbb{F}(1,2)$ is nearly K\"ahler with reduced characteristic holonomy,
then the forms $\tilde{\Omega} =  \tilde{\Omega}_1  +  
\tilde{\Omega}_2$ and $\tilde{\Sigma}$ exist and  $d(\tilde{\Omega}_1 - 
2 \tilde{\Omega}_2) = 0$ holds.
The equation
\bdm
d e_7 \ = \ (2a \, + \, c) \, ( \tilde{\Omega}_1 \, - \, 2 \,
\tilde{\Omega}_2 )
\edm
defines---under the obvious integral condition for the cohomology
class $(2a \, + \, c) \, ( \tilde{\Omega}_1 \, - \, 2 \,
\tilde{\Omega}_2 )$---a $\mathrm{S}^1$-principal
bundle $\pi : M^7 \longrightarrow
\tilde{X}^6$ together with a connection. Finally, 
$M^7$ admits a $\G_2$-structure
\bdm
\varphi \ = \ \pi^*(\tilde{\Sigma}) \, + \, \pi^*(\tilde{\Omega}_1 \, + \, 
\tilde{\Omega}_2) \wedge e_7
\edm
with parallel characteristic torsion form
\bdm
\T^{\mathrm{c}} \ = \ \frac{7}{4} \, c \, \pi^*(\tilde{\Sigma}) \, + \, 
(2 a \, + \, c) \, \pi^*(\tilde{\Omega}_1 \, - \, 2 \, 
\tilde{\Omega}_2) \wedge e_7 \ .
\edm
All together we classified this type of regular $\G_2$-manifolds.
\begin{thm}
Let $(M^7, g, \varphi)$ be a complete, cocalibrated $\G_2$-manifold such that
\bdm
\nabla^{\mathrm{c}} \T^{\mathrm{c}} \ = \ 0  , \quad 
\hol(\nabla^{\mathrm{c}}) = \un(2)
\edm
and suppose that $\T^{\mathrm{c}}$ acts with opposite eigenvalues
$\pm 7 \, c \neq 0$ on the $\nabla^{\mathrm{c}}$-parallel spinors
$\Psi_1 \, , \, \Psi_2$. Moreover, suppose that $M^7$ is regular. Then
$M^7$ is a principal $\mathrm{S}^1$-bundle and a Riemannian submersion
over  the projective space $\mathbb{CP}^3$ or 
the flag manifold $\mathbb{F}(1,2)$ equipped with their standard nearly
K\"ahler structure coming from the twistor construction. The Chern class
of the fibration $\pi : M^7 \longrightarrow  \mathbb{CP}^3 \, , 
\mathbb{F}(1,2)$ is proportional to the K\"ahler form. Conversely,
any of these fibrations admits a $\G_2$-structure with parallel 
characteristic torsion and characteristic holonomy contained in $\un(2)$. 
\end{thm}

\noindent
Consider now the case that $a + b = c$. The torsion form is given by
\bdm
\T^{\mathrm{c}} \ = \ (2 a \, + \, c) \, \Omega_1 \wedge e_7 \, + 
\, ( 5 c - 4 a) \, \Omega_2 \wedge e_7
\edm 
and it acts on the spinors $\Psi_1$ and $\Psi_2$ with 
eigenvalue $- \, 7 c$. Again, we may simplify the formulas
for the derivatives,
\begin{eqnarray*}
de_7 &=& (2 a \, + c) \, \Omega_1 \, + \, (5 c - 4 a)\, \Omega_2 \, , \quad
d\Omega_1 \ = \ d \Omega_2  \ = \ 0 \, , \\
d* \Sigma &=& 0 \, , \quad  d(e_7 \haken * \Sigma) \ = \ 
7 c \, \Sigma \wedge e_7 \, , \quad
d \Sigma \ = \ 7 c \, (* \Sigma)  .
\end{eqnarray*}
In particular, the Lie derivatives are given by 
\bdm
\mathcal{L}_{e_7} \Omega_1 \ = \ 0 \, , \quad
\mathcal{L}_{e_7} \Omega_2 \ = \ 0 \, , \quad
\mathcal{L}_{e_7} \Sigma \ = \ 7 c \, (e_7 \haken * \Sigma) \, , \quad
\mathcal{L}_{e_7} (e_7 \haken * \Sigma) \ = \ - \, 7 c \, \Sigma  .
\edm
The tangent bundle of $M^7$ splits into two complex bundles and one
real bundle,
\bdm
TM^7 \ = \ E_1 \, \oplus \, E_2 \, \oplus \, \R \cdot e_7  ,
\edm
where $E_1$ is spanned by $\{e_1,  e_2 , e_3, e_4\}$ and
$E_2$ is spanned by $\{e_5,e_6\}$. The characteristic connection
preserves this splitting. Using the formula for the characteristic torsion,
we see that the bundles $E_1 \, \oplus \, \R^1 \cdot e_7$ and
 $E_2 \, \oplus \, \R^1 \cdot e_7$ are preserved by the Levi-Civita
connection $\nabla^g$,
\bdm
\nabla^g(E_1 \, \oplus \, \R \cdot e_7) \ \subset\ 
E_1 \, \oplus \, \R^1 \cdot e_7  , \quad 
\nabla^g(E_2 \, \oplus \, \R^1 \cdot e_7) \ \subset\ 
E_2 \, \oplus \, \R^1 \cdot e_7 .
\edm
We compute the Ricci tensor as well as the scalar curvature:
\begin{eqnarray*}
\mathrm{Ric}^{\nabla^{\mathrm{c}}} &=& ( -  4 a^2  +  10 a c  +  
6 c^2) \, \mathrm{Id}_{E_1} \, \oplus \, 
(-  16 a^2  +  12 a c  + 10 c^2)  \mathrm{Id}_{E_2} \ , \\
\mathrm{Ric}^{g} &=& (-  2 a^2  +  12 a c  +  \frac{13}{2}  c^2) 
\mathrm{Id}_{E_1} \, \oplus \, (-  8 a^2  -  8 a c + 
\frac{45}{2}  c^2)  \mathrm{Id}_{E_2}  \\
&\oplus& 
(12 a^2  -  16 ac  +  \frac{27}{2}  c^2) \cdot e_7  , \\
\mathrm{Scal}^g &=& -  12 a^2  +  16 a c  + 
\frac{169}{2}  c^2  , \quad
\mathrm{Scal}^{\nabla^{\mathrm{c}}} \ = \ -  48 a^2  +  64 a c  +  
44 \, c^2  .
\end{eqnarray*}
For a regular structure, the orbit space  $\pi : M^7 \longrightarrow
\tilde{X}^6$ admits a Riemannian metric $\tilde{g}$, two closed forms
$\tilde{\Omega}_1 , \, \tilde{\Omega}_2$ and a splitting of the tangent bundle,
\bdm
T\tilde{X}^6 \ = \  \tilde{E}_1 \, \oplus \, \tilde{E}_2  .
\edm
The Levi-Civita connection of $\tilde{X}^6$ preserves this
splitting. Consequently, the universal covering of 
$\tilde{X}^6$ splits into the product of two
K\"ahler-Einstein 
manifolds $(\tilde{Y}_1, \tilde{g}, \tilde{\Omega}_1)$ and 
 $(\tilde{Y}_2, \tilde{g}, \tilde{\Omega}_2)$. The fibers of the 
Riemannian submersion  $\pi : M^7 \longrightarrow
\tilde{X}^6$ are totally geodesic and the O'Neill tensor is given
by
\bdm
\mathrm{A}_X Y \ = \ \frac{1}{2} \, \mathrm{pr}_{E_2} [X \, , \, Y] \ = \ 
- \, \frac{1}{2} \, \T^{\mathrm{c}}(X , Y , *)  .
\edm
We apply  formula $(9.36c)$ of \cite{Besse} and  compute the Ricci
tensor of $\tilde{X}^6$
\bdm
\tilde{\mathrm{Ric}}^{\tilde{g}} \ = \ 
7c \, (2 a + c) \,
\mathrm{Id}_{\tilde{E}_1} \, \oplus \, 7c \, ( 5 c - 4 a) \, 
\mathrm{Id}_{\tilde{E}_2} .
\edm
Denote by $\tilde{S}_i$ the scalar curvature of $\tilde{Y}_i$ for $i=1,2$.
The sum $\tilde{S} = \tilde{S}_1 \, + \,\tilde{S}_2$ 
is the scalar curvature of  $\tilde{X}^6$, and we have
\bdm
\tilde{S}_1 \ = \ 28 c \, ( 2 a \, + \,  c) \, , \quad
\tilde{S}_2 \ = \ 14 c \, ( 5 c \, - \, 4 a ) \, , \quad
\tilde{S} \ = \ 7 \cdot 14 \cdot c^2 \ > \ 0  .
\edm
By inverting these expression, we may express the parameters $a,c$  by the 
scalar curvatures,
\bdm
a \ = \ \frac{ 5 \, \tilde{S}_1  -  2  \tilde{S}_2}{28 \sqrt{2 \tilde{S} }}
 , \quad c \ = \ \frac{\sqrt{\tilde{S}}}{7 \sqrt{2}}  .  
\edm
The differential $de_7$ is given by
\bdm
d e_7 \ = \ (2a \, + \, c)  \Omega_1  + (5 c  -  4 a) \, \Omega_2 \ 
= \ \frac{1}{\sqrt{2  \tilde{S}}}  \big( \frac{1}{2} \, \tilde{S}_1 
\tilde{\Omega}_1  +  \tilde{S}_2  \tilde{\Omega}_2 \big) \ = \ 
\frac{\sqrt{2}}{\sqrt{\tilde{S}}} \, \tilde{\mathrm{Ric}}^{\tilde{g}}  . 
\edm
The form $\Sigma$ is a section in a line bundle. Indeed, let us introduce
the complex-valued form $\sigma :=  \Sigma  +  i \,
(e_7 \haken * \Sigma)$. Then we obtain
\bdm
\mathcal{L}_{e_7} \sigma \ = \ -  7i\,c  \sigma  .
\edm  
Denote by $L$ the length of the closed integral curves of $e_7$.
$\sigma$ is periodic along the integral curve, i.e.~$7 c  L =  2 k\pi$.
Since
\bdm
\sigma(m \cdot e^{2 \pi t i}) \ = \ e^{- \, 7 c L t i} \cdot 
\sigma(m) \ = \  e^{- \, 2 \pi k t i} \cdot \sigma(m) \ ,
\edm
the map $\sigma$ is a section $\Lambda^3_2(\tilde{X}^6) \otimes G_k$, where
$G_k = M^7 \times_{\mathrm{S}^1} \C$ is the associated bundle defined
by the $\mathrm{S}^1$-representation $z \longrightarrow z^k$. The section
$\sigma$ is parallel. The complex-valued $1$-form 
\bdm
\frac{2 \pi}{L} \, e_7 \cdot i \ : \ TM^7 \longrightarrow \R \cdot i
\edm 
is the connection form in the bundle  $\pi : M^7 \longrightarrow
\tilde{X}^6$. The Chern class of this principal bundle is given by
\bdm
c_1^* \ = \ - \, \frac{1}{2 \pi k} \, \tilde{\mathrm{Ric}}^{\tilde{g}} \ = \ 
- \, \frac{1}{k} \, c_1( \tilde{X}^6) \ . 
\edm
Consequently, the curvature of the bundle 
 $\Lambda^3_2(\tilde{X}^6) \otimes G_k$ vanishes automatically. 
Moreover, in the non-simply-connected case the holonomy of the flat bundle  
$\Lambda^3_2(\tilde{X}^6) \otimes G_k$ is trivial. 
\begin{thm}
Let $(M^7, g, \varphi)$ be a complete, cocalibrated $\G_2$-manifold such that
\bdm
\nabla^{\mathrm{c}} \T^{\mathrm{c}} \ = \ 0  , \quad 
\hol(\nabla^{\mathrm{c}}) = \un(2)
\edm
and suppose that $\T^{\mathrm{c}}$ acts with  eigenvalue
$- \, 7 \, c \neq 0$ on the $\nabla^{\mathrm{c}}$-parallel spinors
$\Psi_1  ,  \Psi_2$. Moreover, suppose that $M^7$ is regular. Then
$M^7$ is a principal $\mathrm{S}^1$-bundle and a Riemannian submersion
over  a K\"ahler manifold $\tilde{X}^6$. This manifold
has the following properties:
\begin{enumerate}
\item The universal covering of $\tilde{X}^6$ splits into a $4$-dimensional
K\"ahler-Einstein manifold and a $2$-dimensional surface with constant 
curvature.
\item The scalar curvature  $\tilde{S} = \tilde{S}_1  + \tilde{S}_2 > 0$
is positive.
\item The K\"ahler forms   $\tilde{\Omega}_1$ and  $\tilde{\Omega}_2$
are globally defined on  $\tilde{X}^6$.
\end{enumerate}
The bundle  $\pi : M^7 \longrightarrow
\tilde{X}^6$ is defined by a connection form. Its curvature is proportional 
to the Ricci form of $\tilde{X}^6$.
Finally, the flat bundle $\Lambda^3_2(\tilde{X}^6) \otimes G_k$ admits a
parallel section. Conversely, any $\mathrm{S}^1$-bundle resulting from
this construction admits a cocalibrated $\G_2$-structure such that the
characteristic torsion is parallel and the characteristic holonomy is
contained in $\un(2)$.
\end{thm}
\vspace{2mm}

\begin{exa}
Let $\tilde{Y}_1$ be a simply-connected K\"ahler-Einstein manifold
with negative scalar curvature $\tilde{S}_1 = - 1$, 
for example a hypersurface of degree $d \geq 5$ in $\mathbb{CP}^3$.
For the second factor we choose the round sphere normalized by the 
condition   $\tilde{S}_2 = + 2$. Then the product $\tilde{X}^6 = 
\tilde{Y}_1 \times \tilde{Y}_2$ is simply-connected and the
$\mathrm{S}^1$-bundle defined by the Ricci form
admits a cocalibrated $\G_2$-structure with parallel torsion. Since
the  product $\tilde{X}^6 = 
\tilde{Y}_1 \times \tilde{Y}_2$ is simply-connected, the flat bundle
 $\Lambda^3_2(\tilde{X}^6) \otimes G_1$ admits a  parallel
section $\sigma$. 
\end{exa}
\vspace{2mm}

\noindent
Finally we study the case $c = 0$. The $\G_2$-manifold is of pure type
$\mathcal{W}_3$. In this case the $3$-form $\Sigma$ projects onto
$\tilde{X}^6$
and defines a parallel form in  $\Lambda^3_2(\tilde{X}^6)$. On the other hand,
the curvature of the bundle $\pi : M^7 \rightarrow
\tilde{X}^6$ is proportional to $\tilde{\Omega}_1- 2\tilde{\Omega}_2$.
\begin{thm}
Let $(M^7, g, \varphi)$ be a complete $\G_2$-manifold of pure type
$\mathcal{W}_3$ such that
\bdm
\nabla^{\mathrm{c}} \T^{\mathrm{c}} \ = \ 0  , \quad
 \T^{\mathrm{c}} \ \neq 0  , \quad 
\hol(\nabla^{\mathrm{c}}) = \un(2)  .
\edm
Moreover, suppose that $M^7$ is regular. Then
$M^7$ is a principal $\mathrm{S}^1$-bundle and a Riemannian submersion
over  a Ricci-flat K\"ahler manifold $\tilde{X}^6$. This manifold
has the following properties:
\begin{enumerate}
\item The universal covering of $\tilde{X}^6$ splits into a $4$-dimensional
Ricci-flat K\"ahler manifold and the $2$-dimensional flat space $\R^2$.
\item The K\"ahler forms   $\tilde{\Omega}_1$ and  $\tilde{\Omega}_2$
are globally defined on  $\tilde{X}^6$.
\item There exists a parallel form $\Sigma \in \Lambda^3_2(\tilde{X}^6)$.
\end{enumerate}
The bundle  $\pi : M^7 \longrightarrow
\tilde{X}^6$ is defined by a connection form. Its curvature is proportional
to the form 
\bdm
\tilde{\Omega}_1 \, - \, 2 \, \tilde{\Omega}_2 .
\edm
Conversely, any $\mathrm{S}^1$-bundle resulting from
this construction admits a $\G_2$-structure of pure type 
$\mathcal{W}_3$ such that the
characteristic torsion is parallel and the characteristic holonomy is
contained in $\un(2)$.
\end{thm}
\vspace{2mm}

\begin{exa}
Consider a $K3$-surface and denote by $\tilde{\Omega}_1$ its K\"ahler
form. Then there exist two parallel forms $\eta_1  ,  \eta_2$ 
in $\Lambda^2_+(K3)$ being orthogonal to  $\tilde{\Omega}_1$. Let $e_5$
and $e_6$ be a parallel frame on the torus $T^2$. The product
$\tilde{X}^6 = K3 \times T^2$ satisfies the conditions of the latter
Theorem. Indeed, we can construct the following parallel form
\bdm
\Sigma \ = \ \eta_1 \wedge e_5  +  \eta_2 \wedge e_6  .
\edm
Moreover, the cohomology class of $\tilde{\Omega}_1  -  2 \tilde{\Omega}_2$ 
has to be
proportional to an integral class. This implies the condition that
$\tilde{\Omega}_1/\mathrm{vol}(T^2) \in  \mathrm{H}^2 ( K3  ; \Q)$ 
is a rational cohomology class.
\end{exa}
%
\section{$\G_2$-manifolds with parallel torsion and 
$\hol(\nabla^{\mathrm{c}}) = 
\su(2)$}\noindent
We briefly discuss the structure of simply-connected, complete, cocalibrated
$\G_2$-manifolds with parallel characteristic torsion and 
$\hol(\nabla^{\mathrm{c}}) = \su(2)$. The tangent bundle 
splits into the sum of two bundles preserved by the
characteristic connection,
\bdm
TM^7 \ = \ E_1 \, \oplus \, E_2  .
\edm
In our notation, the three-dimensional  subbundle $E_2$ is spanned
by $\{e_5 \, , \, e_6 \, , \, e_7\}$. Moreover, the following forms
are globally defined and $\nabla^{\mathrm{c}}$-parallel,
\bdm
e_5 \, , \quad e_6 \, , \quad e_7 \, , \quad 
\Omega_1 \ := \, e_{12} \, + \, e_{34} \, , \quad
\Omega_2 \ := \, e_{14} \, + \, e_{23} \, , \quad
\Omega_3 \ := \, e_{13} \, - \, e_{24} \ .
\edm
The $\G_2$-structure is given by the formula
\bdm
\varphi \ = \ \Omega_1 \wedge e_7  -  \Omega_2 \wedge e_6  + 
\Omega_3 \wedge e_5 +  e_{567}  .
\edm
Basically there are three algebraic types of
torsion forms. If
$\T^{\mathrm{c}} = e_{567}$, then $M^7$ splits into the product
$M^7 = Y^4 \times \mathrm{S}^3$ of the sphere $\mathrm{S}^3$ by 
a simply-connected, complete,
Ricci-flat and anti-selfdual manifold $Y^4$. The forms
$\Omega_1  ,  \Omega_2  ,  \Omega_3$ are the parallel forms
in $\Lambda^2_+(Y^4)$. Conversely, any product of that type admits
a cocalibrated $\G_2$-structure with holonomy $\hol(\nabla^{\mathrm{c}}) = 
\su(2)$. If $\T^{\mathrm{c}} =  \Omega_3 \wedge e_5 $, 
then $e_6  ,  e_7$ are
$\nabla^g$-parallel. The manifold splits into $M^7 = Y^5 \times
\R^2$. Moreover, we obtain
\bdm
d e_5 \ = \  \Omega_3  , \quad  
(d e_5)^2 \wedge e_5 \ \neq \ 0  , \quad
\T^{\mathrm{c}} \ = \ de_5 \wedge e_5.
\edm
The tuple $(Y^5, g , e_5, \Omega_3)$ is homothetic to a Sasakian manifold
with characteristic torsion $\T^{\mathrm{c}}$ and holonomy
$\su(2)$. These spaces have been described completely in
\cite {FriedrichIvanov} , Theorem 7.3. and Example 7.4. They are 
$\eta$-Einstein Sasakian manifolds with Ricci tensor $\mathrm{Ric}^g = (6 , 6,
6 , 6, 4)$. Again, we can reconstruct the $\G_2$-structure of $Y^5 \times
\R^2$ out of the $\eta$-Einstein Sasakian structure of $Y^5$. The third
possibility for the torsion form is
\bdm
\T^{\mathrm{c}} \ = \  \Omega_1 \wedge e_7  +  \Omega_2 \wedge e_6  - 
\Omega_3 \wedge e_5  - 2   e_{567}  .
\edm
A computation of the Ricci tensor yields the following result:
\bdm
\mathrm{Ric}^{\nabla^{\mathrm{c}}} \ = \ 3 \,
\mathrm{Id}_{E_1} \, \oplus 0 \, \mathrm{Id}_{E_2}  , \quad
\mathrm{Ric}^{g} \ = \ \frac{9}{2} \,
\mathrm{Id}_{E_1} \, \oplus 3  \mathrm{Id}_{E_2}  .
\edm
Since $e_5  ,  e_6  ,  e_7$ are $\nabla^{\mathrm{c}}$-parallel, the 
Killing vector fields define a locally free isometric
action of the group $\SU(2)$ on $M^7$. We identify the Lie algebra
of $\SU(2)$ with these Killing vector fields.  
In the regular case the $\G_2$-manifold
is a principal $\SU(2)$-bundle $M^7 \rightarrow Y^4$ over a smooth Riemannian
four-manifold. The vector valued $1$-form $Z \, : \, TM^7 \rightarrow \su(2)$ 
defined by
\bdm
Z \ := \ e_5 \, \otimes \, e_5 \, + \, e_6 \, \otimes \, e_6 \, + \, 
e_7 \, \otimes \, e_7
\edm
is a connection form. The formulas
\bdm
d e_5 \ = \ - \, \Omega_3 \, - \, 2 \, e_{67} \, , \quad
d e_6 \ = \ \Omega_2 \, + \, 2 \, e_{57} \, , \quad
d e_7 \ = \ \Omega_1 \, - \, 2 \, e_{56} \, , 
\edm
express the curvature of the connection $Z$,
\bdm
\Omega_Z \ = \ \Omega_1  \otimes e_7 + \Omega_2 \, \otimes \, e_6  - 
\Omega_3 \otimes  e_5  . 
\edm
Consequently, $M^7$ is a $\SU(2)$-instanton bundle over $Y^4$ and the selfdual
curvature is parallel.
%
%
\section{$\G_2$-manifolds with parallel torsion and 
$\hol(\nabla^{\mathrm{c}}) = 
\so(3)$}\noindent
A $\G_2$-manifold with characteristic holonomy $\hol(\nabla^{\mathrm{c}}) = 
\so(3)$ admits two $\nabla^{\mathrm{c}}$-parallel spinor fields $\Psi_1 \, ,
\, \Psi_2$, a $\nabla^{\mathrm{c}}$-parallel vector field
$e_7$ and a $\nabla^{\mathrm{c}}$-parallel $2$-form $e_{12} + e_{34} +
e_{56}$.  Moreover,
the representation splits into $\R^7 = \R^3 \oplus \R^3 \oplus \R^1$.
Consequently,
the condition $\nabla^{\mathrm{c}} \T^{\mathrm{c}} = 0$ implies that the
manifold is naturally reductive (see \cite{CleytonSwann}). The manifold is a
homogeneous space $\G/\SO(3)$, where the Lie algebra $\g$ 
of the $10$-dimensional
group $\G$ is completely fixed by the torsion and the curvature of the 
characteristic connection. Indeed, we have $\g = \so(3) \oplus \R^7$ and
the bracket is given by the formula
\bdm
\big[ A \, + \, X \, , \, B \, + \, Y \big] \ = \ 
\big( [A \, , \, B] \, - \, \mathrm{R}^{\mathrm{c}}(X , Y) \big) \, + \, 
\big( A \cdot Y \, - \, B \cdot X \, - \, \T^{\mathrm{c}}(X , Y) \big)  .
\edm
The family of admissible torsion forms depends on three parameters,
\begin{eqnarray*}
\T_{a,b,c} &=& d \cdot \varphi \, + \,  
a\, (- 2 \, e_{123} \, + \, e_{136} \, - \, e_{145} \, + \,
e_{235} \, + \, e_{246} \, + \, 2 \, e_{356}) \\ 
&+& \, b \, (- 2 \, e_{124} 
\, - \, e_{135} \, - \, e_{146} \, + \,
e_{236} \, - \, e_{245} \, + \, 2 \, e_{456}) \\ 
&+& \, c \, (4 \, e_{127} 
\, - \, 3 \, e_{135} \, + \, 3\, e_{146} \, + \,
3 \, e_{236} \, + \, 3 \, e_{245} \, + \, 4 \, e_{347} \, + \, 4 \, e_{567})
\ ,
\end{eqnarray*} 
where either $d = 3 \, c$ or $d = - \, 4 \, c$. If $d = 3c$,  the torsion
form acts on $\Psi_1, \Psi_2$ with the same eigenvalue. For $d = - 4c$ it acts
on these spinors with opposite eigenvalues. The curvature operator
$\mathrm{R}^{\mathrm{c}} : \Lambda^2(\R^7) \rightarrow \so(3)$ is invariant.
This gives a five-dimensional space parameterized by $\{x,y,z,u,v\}$,
\begin{eqnarray*}
\mathrm{R}^{\mathrm{c}} &=& \mathrm{R}^{\mathrm{c}}_1 \otimes S_1 
\, + \, \mathrm{R}^{\mathrm{c}}_2 \otimes S_2 \, + \, 
 \mathrm{R}^{\mathrm{c}}_3 \otimes S_3 \, , \\
\mathrm{R}^{\mathrm{c}}_1 &=&  (x - y) (e_{12} - e_{56})  +  
(x + y) (e_{16} + e_{25}) +  2 \, z (e_{26} - e_{15})  +  2 \, u \,
e_{37} +  2 \, v \, e_{47} \, , \\
\mathrm{R}^{\mathrm{c}}_2 &=&  x(e_{13} + e_{35}) + 
z (e_{14} + e_{23} + e_{36} + e_{45})
+ v (e_{17} + e_{57}) - y (e_{24} + e_{46}) \\
&-& u (e_{27} + e_{67}) \, , \\
\mathrm{R}^{\mathrm{c}}_3 &=& z (e_{13} - e_{24} - e_{35} + e_{46}) 
+ y (e_{45} - e_{14} ) + u (e_{57} - e_{17}) + x (e_{36} - e_{23}) \\
&+& v (e_{67} - e_{27}) \ . 
\end{eqnarray*}
Since the characteristic torsion is parallel, the curvature tensor
$\mathrm{R}^{\mathrm{c}}(X,Y,U,V)$ is symmetric with respect to
the pairs $(X,Y)$ and $(U,V)$, see \cite{FriedrichIvanov}. This implies 
directly that $u = v = z = 0$ and $y = -  x$. The curvature
operator is proportional to the projection $\mathrm{pr} : \Lambda^2(\R^7)
\rightarrow \so(3)$ onto the Lie subalgebra,
\bdm
\mathrm{R}^{\mathrm{c}} \ = \ x  \big(2 S_1 \otimes S_1 
 +   S_2 \otimes S_2  +  S_3 \otimes S_3 \big)  .
\edm
The pair $(\T^{\mathrm{c}} , \mathrm{R}^{\mathrm{c}})$ has to satisfy the
Bianchi idendity. In particular, $\T^{\mathrm{c}} \cdot \T^{\mathrm{c}} 
+ \mathrm{R}^{\mathrm{c}}$ is a scalar in the Clifford algebra.
This equation has two solutions, namely:
\begin{eqnarray*}
d &=& -  4  c  , \quad x \ = \ a^2  +  b^2  -  49  c^2  , \\
d &=& 3  c  , \quad a \ = \ b \ = \ 0  ,  
\quad 2  x \ = \ -  49  c^2  .
\end{eqnarray*}
If $d = 3  c \neq 0$, then the parameters $a = b = 0$
vanish. The torsion form depends only on the parameter $c$ and can
be written in the simpler form
\bdm
\T^{\mathrm{c}} \ = \ 7 \, c \, (e_{12} \, + \, e_{34} \, + \, e_{56}) \wedge
e_7 \, , \quad 
\mathrm{R}^{\mathrm{c}} \ = \ - \, \frac{49}{2} \, c^2 \, \big(2 \, 
S_1 \otimes S_1 
\, + \,  S_2 \otimes S_2 \, + \, S_3 \otimes S_3 \big) \ .
\edm
A computation of the Ricci tensor yields the following result
\bdm
\mathrm{Ric}^g(e_i) \ = \ \frac{5}{2} \, 49 \, c^2 \, e_i \quad
\mbox{for} \quad i \ = \ 1 \, , \ldots , \, 6 \, , \quad
\mathrm{Ric}^g(e_7) \ = \ \frac{3}{2} \, 49 \, c^2 \, e_7 \ .
\edm
In particular, the Riemannian Ricci tensor is positive definite and,
consequently, $M^7$ is compact. One easily identifies the group
$\G$: It is the group $\G = \SO(5)$ with the standard
embedding of $\SO(3)$. The corresponding naturally reductive space is the
Stiefel manifold $M^7 = \SO(5)/\SO(3)$.
\vspace{2mm}

\begin{thm}
A simply-connected, complete, cocalibrated $\G_2$-manifold with characteristic
holonomy $\hol(\nabla^{\mathrm{c}}) = \so(3)$ such that $\T^{\mathrm{c}}$
acts with the same eigenvalue on the parallel spinors $\Psi_1 , \Psi_2$
is isometric to the Stiefel manifold $\SO(5)/\SO(3)$. 
The metric
is a Riemannian submersion over the Grassmanian manifold $\G_{5,2}$.
\end{thm} 
\vspace{2mm}

\noindent
If $d = - \, 4 \, c$,  the torsion form does not depend on
$e_7$, $e_7 \haken \T^{\mathrm{c}} = 0$. Since $e_7$ is 
$\nabla^{\mathrm{c}}$-parallel, the vector field is parallel
with respect to the Levi-Civita connection, too. A complete, simply-connected
$\G_2$-manifold of that type splits into the Riemannian product $Y^6 \times
\R^1$, where $Y^6$ is an almost Hermitian manifold of Gray-Hervella-type 
$ \mathcal{W}_1 \oplus \mathcal{W}_3$ with
characteristic holonomy $\so(3) \subset \su(3)$. $Y^4$ is completely
defined by the torsion form and this family of almost Hermitian 
$6$-manifolds has
been studied in \cite{AFS}, Theorem 4.6 as well as in \cite{Sch}.
Let us summarize the result.
\vspace{2mm}

\begin{thm}
A simply-connected, complete, cocalibrated $\G_2$-manifold with characteristic
holonomy $\hol(\nabla^{\mathrm{c}}) = \so(3)$ such that $\T^{\mathrm{c}}$
acts with opposite eigenvalues on the parallel spinors $\Psi_1 , \Psi_2$
splits into  
the Riemannian product $Y^6 \times
\R^1$, where $Y^6$ is an almost Hermitian manifold of Gray-Hervella-type 
$ \mathcal{W}_1 \oplus \mathcal{W}_3$ with
characteristic holonomy $\so(3) \subset \su(3)$. 
\end{thm} 
\section{$\G_2$-manifolds with parallel torsion and 
$\hol(\nabla^{\mathrm{c}}) = 
\so_{ir}(3)$}\noindent
Since $(\Lambda^3_{27})_{\so_{ir}(3)}= 0$ is trivial, any cocalibrated
$\G_2$-manifold with parallel characteristic torsion and 
$\hol(\nabla^{\mathrm{c}}) = \so_{ir}(3)$ is nearly parallel. Moreover, the
curvature tensor $\mathrm{R}^{\mathrm{c}}$ is $\nabla^{\mathrm{c}}$-parallel, 
see \cite{CleytonSwann}. There exists only one $\so_{ir}(3)$-invariant
curvature operator $\mathrm{R}^{\mathrm{c}} : \Lambda^2(\R^7) \longrightarrow
\so_{ir}(3)$, namely the projection onto the subalgebra $\so_{ir}(3) 
\subset \so(7) = \Lambda^2(\R^7)$. Consequently, the characteristic torsion
and the curvature operator are uniquely defined. On the other side, 
consider the embedding of $\SO(3)$ into $\SO(5)$
given by the $5$-dimensional irreducible $\SO(3)$-representation. Then
the naturally reductive space $\SO(5)/\SO_{ir}(3)$ admits a nearly parallel
$\G_2$-structure, see \cite{FKMS}. Finally 
we obtain a complete classification in this case.
\begin{thm} \label{irrso(3)}
A complete, simply-connected and cocalibrated $\G_2$-manifold with parallel
characteristic torsion and  $\hol(\nabla^{\mathrm{c}}) = 
\so_{ir}(3)$ is isometric to $\SO(5)/\SO_{ir}(3)$.
\end{thm}
%
%
\section{$\G_2$-manifolds with parallel torsion and 
$\hol(\nabla^{\mathrm{c}}) = 
\su_{c}(2)$}\noindent
In this section we prove the following uniqueness result.
\begin{thm} \label{Nonexistence}
There exists a unique simply-connected, complete, 
cocalibrated $\G_2$-manifold with
\bdm
\nabla^{\mathrm{c}} \T^{\mathrm{c}} \ = \ 0  , \quad
\hol(\nabla^{\mathrm{c}}) \ = \ \su_c(2) .
\edm 
The manifold is homogeneous naturally reductive.
\end{thm}
\begin{proof}
The $\su_c(2)$-representation $\R^7 = \C^2 \oplus \R^3$ splits into the
sum of two irreducible representations. Then any $\G_2$-structure
with parallel torsion and holonomy $\su_c(2)$ is naturally reductive,
$\nabla^{\mathrm{c}} \T^{\mathrm{c}} =  0  , 
\nabla^{\mathrm{c}} \mathrm{R}^{\mathrm{c}} =  0$, see \cite{CleytonSwann}.
The torsion forms belong to $\mathrm{Tor}_{\su_c(2)}$ and are
parameterized by two parameters,
$\T^{\mathrm{c}} =  a \cdot \varphi \, + \, b \cdot e_{567}$.
The curvature operator $\mathrm{R}^{\mathrm{c}} : \Lambda^2(\R^7) 
\longrightarrow \su_c(2) \subset \Lambda^2(\R^7)$ 
is symmetric and   $\su_c(2)$-invariant. Consequently, the curvature
operator is proportional to the projection onto the Lie subalgebra
$\su_c(2)$,
\begin{eqnarray*}
\mathrm{R}^{\mathrm{c}} \ = \ p  ( Q_1 \otimes Q_1  +   
Q_2 \otimes Q_2 + Q_3 \otimes Q_3 )  .
\end{eqnarray*}
$Q_1, \, Q_2$ and $Q_3$ denote the basis of the Lie algebra
$\su_c(2)$ introduced before. The pair $(\T^{\mathrm{c}} \, , \, 
\mathrm{R}^{\mathrm{c}})$ satisfies the Bianchi identity if and only if
$(\T^{\mathrm{c}})^2 + \mathrm{R}^{\mathrm{c}}$
is a scalar in the Clifford algebra, see \cite{Kostant}. There are two
solutions of this algebraic equation,
\bdm
a \ = \ 0  , \ p \ = \ 0  , \quad \mbox{and} \quad p \ = \ a^2 , 
5  a  +  b \ = \ 0 .
\edm
The case $p = 0  , a = 0$ defines a flat structure, 
$\hol(\nabla^{\mathrm{c}}) = 0$. The second case yields a unique
naturally reductive, cocalibrated $\G_2$-manifold with parallel
torsion and holonomy $\hol(\nabla^{\mathrm{c}}) = \su_c(2)$  
The Lie algebra $\g$ of the $10$-dimensional automorphism group 
is given by
$\g = \su_c(2) \oplus \R^7$ with the bracket
\bdm
\big[ A \, + \, X \, , \, B \, + \, Y \big] \ = \ 
\big( [A \, , \, B] \, - \, \mathrm{R}^{\mathrm{c}}(X , Y) \big) \, + \, 
\big( A \cdot Y \, - \, B \cdot X \, - \, \T^{\mathrm{c}}(X , Y) \big) \, .
\edm
It turns out that $\g$ is perfect, $[\g \, , \, \g] = \g$. 
The adjoint representation $\mathrm{ad} : 
\g \rightarrow \mathfrak{gl}(\g)$ is a faithful representation of $\g$. 
The Lie
algebra $\g$ contains a $7$-dimensional nipotent radical $\mathfrak{r}$. It
is generated by $ \mathfrak{r} =
\mathrm{Lin}(e_1,e_2,e_3,e_4,e_5-Q_2,e_6+Q_1,e_7+Q_3)$. Moreover, 
$[\mathfrak{r} \, , \,\mathfrak{r}] =
\mathrm{Lin}(e_5-Q_2,e_6+Q_1,e_7+Q_3)$ is three-dimensional and abelian.
$\g/\mathfrak{r} = \su_c(2)$ is isomorphic to the holonomy
algebra. 
\end{proof}
\noindent
The case of $\hol(\nabla^{\mathrm{c}}) = \R^1 \oplus \su_c(2)$ is similar.
The admissible torsion forms are again $\T^{\mathrm{c}} =  a \cdot \varphi \,
+ \, b \cdot e_{567}$ and the 
$(\R^1 \oplus \su_c(2))$-invariant operators  $\mathrm{R}^{\mathrm{c}} : 
\Lambda^2(\R^7) \longrightarrow \R^1 \oplus \su_c(2)$ 
are parameterized by three parameters, 
\begin{eqnarray*}
\mathrm{R}^{\mathrm{c}} &=& p \cdot \big( e_{34} \otimes Q_3 \, - \, e_{14} 
\otimes Q_1 \, - \, e_{23} \otimes Q_1 \, - \, e_{13} \otimes Q_2 \, + \, 
e_{24} \otimes Q_2 \, + \, e_{12} \otimes Q_3 \big) \\ 
&+& q \cdot \big( e_{56} \otimes Q_3 \, - \, e_{57} \otimes Q_1 \, - \,
e_{67} \otimes Q_2 \big) \, + \, r \cdot \big(e_{13} \, + \, e_{24} \big)
\otimes P_1 \ .
\end{eqnarray*}
Since $\mathrm{R}^{\mathrm{c}}$ is symmetric, we have $ q = - \, 2 \, p$. 
The curvature operator simplifies,
\bdm
\mathrm{R}^{\mathrm{c}} \ = \ p \, ( Q_1 \otimes Q_1 \, + \,  
Q_2 \otimes Q_2 \, + \,Q_3 \otimes Q_3 ) \, + \, r \, P_1 \otimes P_1 \ .
\edm
$(\T^{\mathrm{c}})^2 \, + \,\mathrm{R}^{\mathrm{c}}$
is a scalar in the Clifford algebra if and only if the following relations
between the parameters hold,
\bdm
p \ = \ - \, \frac{1}{2} \, a \, (3a \, + \, b) \, , \quad
r \ = \ \frac{3}{2} \, a \, (5a \, + \, b) \ .
\edm
If $5a \, + \, b \neq 0$ and  $3a \, + \, b \neq 0$, then the holonomy
of the characteristic connection is the full Lie algebra 
$\R^1 \oplus \su_c(2)$. All together we classified  simply-connected,
complete, cocalibrated $\G_2$-manifolds with parallel characteristic torsion
and holonomy $\R^1 \oplus \su_c(2)$. The spaces are naturally reductive.
Up to a scaling, the family depends on one parameter.
\begin{NB}
If $b = 0$, then $M^7 = N(1,1)$ is a nearly parallel $\G_2$-manifold with an
$11$-dimensional automorphism group. The automorphism group is isomorphic
to $\SU(3) \times \SU(2)$ and the space appears
in the classification
of all nearly parallel $\G_2$-manifolds with a large automorphism group,
see \cite{FKMS}.
\end{NB}
%
\section{$\G_2$-manifolds with parallel torsion and 
$\hol(\nabla^{\mathrm{c}}) = 
\su(2) \oplus \su_c(2)$}\noindent
\begin{exa}
Starting with a $3$-Sasakian manifold and rescaling again its
metric along the three-dimensional bundle
spanned by  $e_5 , e_6 , e_7$, one obtains a family $(M^7 , g_s , \varphi_s)$ 
of cocalibrated $\G_2$-manifold such that
\bdm
d \, *_s \, \varphi_s \ = \ 0  , \quad \T_s^{\mathrm{c}} \ = \ 
\big( \frac{2}{s} \, - \, 10 s \big) \, e^*_{567}  +  2s  
\varphi_s  , \quad
\nabla^{\mathrm{c}} \T_s^{\mathrm{c}} \ = \ 0  .
\edm
The necessary computations proving these properties are contained in
 \cite{FKMS}, Theorem 5.4.
The characteristic connection preserves the splitting of the tangent
bundle and, consequently, its holonomy is
$\hol(\nabla^{\mathrm{c}}) = 
\su(2) \oplus \su_c(2)$. If $s = 1/\sqrt{5}$, the structure is nearly parallel.
Since $(\T_s^{\mathrm{c}}  ,  \varphi_s) = 4s  +  2/s > 0$, 
these structures are never of pure type $\mathcal{W}_3$.
In the parametrization of the admissible
torsion forms $\T^{\mathrm{c}} = a  \varphi + b  e_{567}$ the family 
realizes the following curve in the $\{a\, , \, b \}$-plane,
\bdm
a \ = \ 2s  , \quad b \ = \ \frac{2}{s}  -  10 s  , \quad
5 a^2 +  ab \ = \ 4 .
\edm
\end{exa}

\noindent
First we investigate naturally reductive $\G_2$-manifolds with holonomy
$\hol(\nabla^{\mathrm{c}}) = \su(2) \oplus \su_c(2)$. The characteristic
torsion and the curvature operator depend on two parameters,
$\mathrm{T}^{\mathrm{c}} = a \cdot \varphi  +  b \cdot e_{567}$ and 
\bdm
\mathrm{R}^{\mathrm{c}} \ = \ p \, ( Q_1 \otimes Q_1 +   
Q_2 \otimes Q_2 \, + \,Q_3 \otimes Q_3 )  +
 r \, ( P_1 \otimes P_1 \, + \,  
P_2 \otimes P_2 \, + \,P_3 \otimes P_3 )  .
\edm
The pair $(\mathrm{T}^{\mathrm{c}} \, , \,\mathrm{R}^{\mathrm{c}})$ satisfies
the Bianchi identity if and only if
\bdm
r \ = \ - \, \frac{a}{2} \, (5 a \, + \, b) \, , \quad p \ = \ 
- \, \frac{a}{2} \, ( 3 a \, + \, b ) .
\edm 
Consequently, we obtain (up to scaling) a one-parameter family of naturally
reductive homogeneous $\G_2$-manifolds with 
$\hol(\nabla^{\mathrm{c}}) = \su(2) \oplus \su_c(2)$.
\vspace{2mm}

\begin{NB}
For $b = 0$, the manifold $M^7$ is nearly parallel and has  a
$13$-dimensional automorphism group, the squashed $7$-sphere.
It appears again in the  classification in  \cite{FKMS}.
\end{NB}
\vspace{2mm}

\noindent
A classification of all cocalibrated 
$\G_2$-manifolds with parallel torsion and holonomy 
$\su(2) \oplus \su_c(2)$ seems to be unaccessible. Nevertheless we can
discuss the geometry of such manifolds and describe some particular cases.
The tangent bundle splits into the sum of two bundles preserved by the
characteristic connection,
\bdm
TM^7 \ = \ E_1  \oplus  E_2  .
\edm
In our notation, the three-dimensional  subbundle $E_2$ is spanned
by $\{e_5  ,  e_6  ,  e_7\}$.
The parallel torsion form of a cocalibrated $\G_2$-manifold depends 
on two parameters,
\bdm
\T \ = \ a \varphi  +  b  e_{567}  .
\edm
\noindent
The  $\nabla^{\mathrm{c}}$-parallel spinor field $\Psi_1$
satisfies the following differential equations
\bdm
\nabla^g_X \Psi_1 \ = \ -  \frac{3a}{4} \, X \cdot \Psi_1 \quad
\text{for} \quad X \in E_1  , \quad
\nabla^g_V \Psi_1 \ = \ -  \frac{3a + b}{4} \, V \cdot \Psi_1 \quad
\text{for} \quad V \in E_2  .
\edm 
\noindent
The Ricci tensor depends only on $\T$ and can be computed explicitly
\begin{eqnarray*}
\mathrm{Ric}^{\nabla^{\mathrm{c}}} &=& (12 a^2 \, + \, 3 a b) \,
\mathrm{Id}_{E_1} \, \oplus (12 a^2 \, + \, 4 a b) \, \mathrm{Id}_{E_2} \ , \\
\mathrm{Ric}^{g} &=& (\frac{27}{2} a^2 \, + \, 3 a b) \,
\mathrm{Id}_{E_1} \, \oplus (13 a^2 \, + \, 4 a b \, + \, \frac{1}{2} \, (a \,
+ \, b)^2) \, \mathrm{Id}_{E_2}.
\end{eqnarray*}
Since $\nabla^{\mathrm{c}}$ preserves the splitting $TM^7 = E_1  \oplus 
E_2$, the algebraic formula for the torsion yields that $E_2$ is an involutive
subbundle. Moreover, the leaves are
totally geodesic. We prove now that every leaf of this distribution is 
a $3$-dimensional Riemannian manifold of constant sectional curvature. The
result is a consequence of the following formulas.
\vspace{2mm}

\begin{lem}
Let $X$ be a vector field in  $E_1$ and $V$ be a vector field
in $E_2$. Then we have
\begin{eqnarray*}
g \big( \nabla^g_V V  ,  \nabla^g_X X \big) &=& 0  , \quad
g \big( \nabla^g_X \nabla^g_V V  , \, X \big) \ = \ 0  , \\
R^g \big(X,V , V , X \big) &=& \frac{1}{4} || \T(X , V , *) ||^2  , \quad
\sum^4_{i = 1} R^g \big(e_i ,  V  , V  ,  e_i \big) \ = \
a^2 \, ||V||^2 .
\end{eqnarray*}
\end{lem} 
\begin{proof}
$\nabla^{\mathrm{c}}$ preserves the splitting $TM^7 = E_1 \oplus E_2$ and
the torsion is totally skew-symmetric. Then we have
\begin{eqnarray*}
0 &=& g \big(\nabla^{\mathrm{c}}_V V \, , \, \nabla^{\mathrm{c}}_X X \big)
\ = \  g \big(\nabla^g_V V \, , \, \nabla^g_X X \big) \ = \ 
- \, g \big(\nabla^g_X \nabla^g_V V \, , \, X \big) \, + \,
X \big( g \big(\nabla^g_V V \, , \, X \big) \big) \\
&=& - \, g \big(\nabla^g_X \nabla^g_V V \, , \, X \big) \, + \,
X \big( g \big(\nabla^{\mathrm{c}}_V V \, , \, X \big) \big) \ = \
- \, g \big(\nabla^g_X \nabla^g_V V \, , \, X \big) \ .
\end{eqnarray*}
We compute the formula for the curvature tensor 
$R^g \big(X,V,V,X \big)$
in a similar way. Alternatively, one can apply formula $(9.29b)$ in 
\cite{Besse}.
\end{proof}
\vspace{2mm}

\noindent
We already know the Ricci tensor of $M^7$ and that
every leaf of the distribution is
totally geodesic. Then the preceding Lemma yields a formula for the Ricci
tensor of the leaves,
\bdm
(12 a^2  +  4 a b  + \frac{1}{2}
(a+ b)^2) \mathrm{Id} \ = \frac{1}{2} (5 a + b )^2 \mathrm{Id} .
\edm
\begin{prop}
Every leaf is a $3$-dimensional Riemannian manifold of constant sectional
curvature
\bdm
\mathrm{k} \ = \ \frac{1}{4} (5 a  +  b)^2 \ \geq \ 0  . 
\edm 
If $M^7$ is complete,  every maximal leaf of the distribution $E_2$
is isometric to a complete space form of non-negative sectional
curvature. 
\end{prop}
\vspace{2mm}

\noindent
If the space of leaves is a smooth manifold,  $\pi : M^7 \rightarrow
Y^4$ is a Riemannian submersion with totally geodesic fibers. The relevant
tensor relating the geometry of $M^7$ and $Y^4$ is
\bdm
\mathrm{A}_X Y \ = \ \frac{1}{2} \, \mathrm{pr}_{E_2} [X  ,  Y] \ = \ 
- \, \frac{1}{2}  \T^{\mathrm{c}}(X , Y , *)  .
\edm
Remark that in our case the tensor $\mathrm{A} : \Lambda^2(E_1) \rightarrow
E_2$ vanishes on $\Lambda^2_{-}(E_1)$. Moreover, if $a \neq 0$, then  
$\mathrm{A} : \Lambda^2_{+}(E_1) \rightarrow E_2$ is an isomorphism.
We apply the formula $(9.36c)$ of \cite{Besse} and we compute the Ricci
tensor of the space of leaves.
\vspace{2mm}

\begin{prop}
The space of leaves $Y^4$ is an Einstein space. Its Ricci tensor is given by
$ 3 a (5 a  +  b)  \, \mathrm{Id}_{TY^4}$.
\end{prop}
\vspace{2mm}

\noindent
If $a = 0$, then the holonomy is reduced to  $\su(2)$. Indeed,
in this case, the Levi-Civita connection
preserves the splitting of the tangent bundle. 
The universal covering of $M^7$ splits into the Riemannian
product. The three-dimensional factor is a sphere $\mathrm{S}^3$. The four-dimensional factor is anti-selfdual and Ricci-flat.
The vector fields $e_5  ,  e_6  ,  e_7$ become
$\nabla^{\mathrm{c}}$-parallel and globally defined. For example, we obtain
\bdm
\nabla^{\mathrm{c}}_X e_5 \ = \ \nabla^{\mathrm{S}^3}_X e_5 \, - \, 
\frac{b}{2} \, e_6 \wedge e_7 (X  , \ *) \ = \ 0  .
\edm
Consequently, the holonomy is $\hol(\nabla^{\mathrm{c}}) = 
\su(2)$. The second interesting case $b = 0$ corresponds to 
nearly parallel $\G_2$-manifolds. According
to Proposition \ref{Nearlyparallel}, Theorem \ref{irrso(3)} 
and Theorem \ref{Nonexistence},
any nearly parallel $\G_2$-manifold different from $\SO(5)/\SO_{ir}(3)$ and 
$N(1,1) = (\SU(3) \times \SU(2))/(\mathrm{S}^1 \times \SU(2))$ has
characteristic holonomy  $\su(2) \oplus \su_c(2)$ or $\g_2$.
At present, only few examples of nearly parallel $\G_2$-manifolds are known;
hence, a complete classification of the case $b=0$ is unaccessible.

    
\end{document}